\documentclass[12pt,a4paper]{article} 
\usepackage[american]{babel}
\usepackage{amsmath}
\usepackage{latexsym}
\usepackage{amssymb}
\usepackage{theorem}
\usepackage{diagrams}
\usepackage{pstcol}
\parskip\medskipamount 
\theoremstyle{change}
\newtheorem{Thm}{Theorem}[section]
\newtheorem{Cor}[Thm]{Corollary}
\newtheorem{Prop}[Thm]{Proposition}

{\theorembodyfont{\rmfamily}
\newtheorem{Num}[Thm]{}
\newtheorem{Rem}[Thm]{Remark}
\newtheorem{Def}[Thm]{Definition}}

\newcommand{\Fin}{{\mathrm{fin}}}

\newcommand{\inv}{^{-1}}

\newcommand{\bbf}{\mathbb} 
\newcommand{\too}{\longrightarrow}
\newcommand{\R}{{}^*\bbf R}

\begin{document} 
 
\title{\bf Asymptotic cones and ultrapowers of Lie groups} 
 
\author{Linus Kramer \& Katrin Tent%
\thanks{Supported by a {\em Heisenberg-Stipendium}}} 
\maketitle

\section{Introduction}

Asymptotic cones of metric spaces were first invented by Gromov. They are metric spaces which capture the 'large-scale structure' of the underlying metric space. Later, van den Dries and Wilkie gave a more general construction of asymptotic cones using ultrapowers. Certain facts about  asymptotic cones, like the completeness of the metric space, now follow rather easily from saturation properties of ultrapowers, and in this survey, we want to present two applications of the van den Dries-Wilkie approach.  Using ultrapowers we obtain an explicit description of the asymptotic cone of a semisimple Lie group. From this description, using semi-algebraic groups and non-standard methods, we can give a short proof of the Margulis Conjecture. In a second application, we use set theory to answer a question of  Gromov.

\section{Definitions}

The intuitive idea behind Gromov's concept of an asymptotic cone was to  look at a given metric space from an 'infinite distance', so that large-scale patterns should become visible. In his original definition this was done by gradually scaling down the metric by factors $1/n$ for $n\in \bbf N$. In the approach by van den Dries and Wilkie, this idea was captured by ultrapowers. Their construction is more general in the sense that the asymptotic cone exists for any metric space, whereas in Gromov's original definition, the asymptotic cone existed only for a rather restricted class of spaces.

Given a divisible ordered abelian group $\Lambda$, we call $(X,d)$ a $\Lambda$-metric space if $d:X\times X\too \Lambda$ satisfies the usual axioms of a metric, i.e. for all $x,y\in X$, $d(x,y)=d(y,x)$, $d(x,y)=0$ if and only if $x=y$, and the triangle inequality holds.

We can now give the definition of asymptotic cone according to van den Dries and Wilkie.

\begin{Def}
\label{cone} Let $(X,d)$ be a metric space and $\mu$ a nonprincipal ultrafilter on $\omega$.
Let $({}^*X,{}^*d)=\prod_\mu(X,d)$ be the ultrapower of $(X,d)$ with respect to $\mu$. Then, clearly, $({}^*X,{}^*d)$ is a $\R$-metric space with $\R=\prod_\mu \bbf R$.
Let $\alpha=(n:\ n\in \omega)_\mu\in \R$ and $o\in X$.
Set 
\[
X^{\langle\alpha\rangle}=\{x\in {}^*X|\ {}^*d(x,o)\leq n\alpha \mbox{ for some } n\in \bbf N\}.
\]
Define $x\approx_\alpha y$ on $X^{\langle\alpha\rangle}$ if ${}^*d(x,y)<\alpha/n$ for all
$n\in\bbf N$.

Then $Cone(X)=X^{\langle\alpha\rangle}/_{\approx_\alpha}$ is called the \emph{asymptotic cone} of $X$. It carries a natural ($\bbf R$-)metric defined by $d(x,y)=std({}^*d(x,y)/\alpha)$ where $std$ denotes the standard part of the element ${}^*d(x,y)/\alpha$.
\end{Def}
Similarly, we can now define the \emph{ultracone} which will also be
useful in our situation. Let
$\R_\Fin=\{t\in {\R} |\ |t|\leq k \mbox{ for some } k\in \bbf N\}\subseteq{}\R$ denote the set of all
finite nonstandard reals.
\begin{Def}
\label{ultracone} With the previous notation, we define $x\approx_\Fin y$ on ${}^*X$ if ${}^*d(x,y)\leq n$ for some $n\in \bbf N$.
The \emph{ultracone} is defined as $UCone(X)={}^*X/_{\approx_\Fin}$ and carries
an $\R/\R_\Fin$-metric.
\end{Def}
Clearly, the constructions are independent of the base point $o\in X$ since in the cone and in the ultracone all elements of $X$ are identified.
While Gromov's definition does not depend on the choice of an ultrafilter, and the asymptotic cone, if it exists, is unique, the definitions given here raise the obvious question to what extend the asymptotic cone and the ultracone depend on the choice of $\mu$. This question will be the focus of Section~\ref{finpres}.

In any case, asymptotic cones only depend on the large-scale structure of the metric space in the following sense:
\begin{Def}
\label{quasiisometries} Let $\Lambda$ be a divisible ordered abelian group, and $(X,d_X)$ and $(Y,d_Y)$ $\Lambda$-metric spaces.
Then $f:X\too Y$ is called a \emph{quasi-isometry} if there are 
constants $L\in\bbf N_{>0}$ and $C\in\Lambda_{\geq 0}$, such that 
\begin{enumerate}
\item $L^{-1}d_X(x,x')-C\leq d_Y(f(x),f(x'))\leq Ld_X(x,x')+C$ for all $x,x'\in X$, and
\item for every $y\in Y$ there is some $x\in X$ such that  $d_Y(f(x),y)\leq C$.
\end{enumerate}
\end{Def}
Special cases of quasi-isometries are, of course, isometries with $L=1,C=0$ and bi-Lipschitz maps with $C=0$. Notice that quasi-isometries need not even be continuous and that any two bounded (e.g. compact) spaces are quasi-isometric.
\begin{Rem}
\label{remark}If $f:X\too Y$ is a quasi-isometry between real metric spaces, then we get obvious induced maps $Cone(f):Cone(X)\too Cone(Y)$, and $UCone(f):UCone(X)\too UCone(Y)$. Note that $Cone(f)$ and
$UCone(f)$ are bi-Lipschitz maps.
In particular, $Cone(X)$ and $Cone(Y)$ are homeomorphic $\mathbb R$-metric spaces.

Note also that two quasi-isometries between  real metric spaces have bounded distance if and only if they agree on the ultracone.
\end{Rem}
We will be mostly interested in the asymptotic cones of two types of metric spaces: finitely generated groups under the word metric, and Lie groups under some left invariant metric, see Section~\ref{Lie}.

\begin{Num}\textbf{The word metric}
\label{wordmetric} In the context of finitely presented groups with polynomial growth, Gromov studied the asymptotic cones of finitely generated groups under the word metric:
Let $\Gamma$ be a group with generating set $S$ and assume
that $S$ is closed under inverses. Then for any $g\in \Gamma\setminus\{1\}$,
we define $|g|$ to be the \emph{word length} of $g$, i.e.
$|g|=min\{n|\ g=x_1\cdots x_n,\ x_i\in S\}$. Put
$|1|=0$. Then one can check that $d(g,h)=|g\inv h|$ for $g,h\in\Gamma$
defines a $\Gamma$-invariant $\bbf Z$-valued
metric, the {\em word metric} on $\Gamma$.
If $S'$ is another finite generating set, and $d'$ its corresponding metric,
then it is easy to see that the identity map is a quasi-isometry from $(\Gamma,d)$
to $(\Gamma,d')$, see \cite{BH} I.8.17(2).
\end{Num}
In general, it is not easy to calculate the asymptotic cone of a given metric space. 
As we noted before, the definition of the asymptotic cone depends on the
ultrafilter $\mu$. However, as the following examples show, in many cases, this ultrafilter does
not really matter.

\begin{Num}\textbf{Examples:}
\begin{enumerate}
\item Let $\Gamma=\bbf Z^n$ be the free abelian group in $n$-generators
equipped with the word metric defined in \ref{wordmetric} 
(which for the canonical generators coincides with the Manhattan Taxi metric,
$d(x,y)=|x_1-y_1|+\cdots+|x_n-y_n|$). It is easy to see (directly from
the definition) that
the asymptotic cone (with respect to the canonical basepoint)
then is the topological group $\bbf R^n$, also with the Manhattan Taxi metric
\cite{Gromov} p.~37.
Hence it is independent of the ultrafilter.

\item Let $\Gamma$ be the free group of rank $n$ with the word metric.
Its asymptotic cone is a complete homogeneous $\bbf R$-tree and again
independent of the ultrafilter. See \cite{Gromov} and
\cite{BH} for more results.

\item Point showed  in \cite{po} that if $\Gamma$ is a finitely generated
nilpotent group, then the asymptotic cone is homeomorphic to $\bbf R^m$,
where $m$ is the sum of the ranks of the descending central series of $G$. 

\end{enumerate}
\end{Num}

\begin{Num}\textbf{Questions:}
Gromov asked the following questions \cite{Gromov} p.~42:
\begin{enumerate}
\item What is the asymptotic cone of a Lie group?
\item Is there a \emph{finitely presented} group, for which different
ultrafilters yield non-homeomorphic asymptotic cones?
\end{enumerate}
\end{Num}
We will show how to determine the asymptotic cone of a Lie group via ultrapowers in the next section. Question (ii) will be seen to be closely related and will be partially answered in Section~\ref{finpres}.
Examples of \emph{finitely generated} groups with
non-homeomorphic asymptotic cones were given by
Thomas and Velickovic \cite{TV}.

\section{Asymptotic cones of Lie groups}\label{Lie}

Let $G$ be a semisimple Lie group of non-compact type. We will also use the notation $G=G(\bbf R)$ to emphasize that we will consider $G$ as the group of $\bbf R$-rational points of some semi-algebraic group defined over the rational numbers. So for any real closed field $R$ we may consider the group $G(R)$. 

Up to quasi-isometry, the group $G(\bbf R)$ carries a unique left-invariant Riemannian metric (see \cite{dlH} IV \S21). We will describe this metric on (the diagonalizable part of) $G$ using the real logarithm. Let $K=K(\bbf R)$ be a (suitably chosen) maximal compact subgroup of $G(\bbf R)$, and let $A=A(\bbf R)$ denote a (semialgebraically connected)
maximal split torus of $G(\bbf R)$, which we may assume to consist of diagonal matrices with positive entries. Then $G(\bbf R)=K(\bbf R)A(\bbf R)K(\bbf R)$, and since $K(\bbf R)$ is compact and we are only concerned with metrics up to quasi-isometry, it suffices to define a $K$-invariant metric on $A(\bbf R)$. 
This can be done in the following way: for $g={\rm diag}(x_1,\ldots x_n)\in A(\bbf R)$ we put $d({\bf 1},g)= (|\log x_1|^2+\ldots +|\log x_n|^2)^{1/2}$. This metric extends to a left-invariant metric $d$ on $G$, i.e.
$d(g,g')=d(\mathbf 1,g^{-1}g')$ for all $g,g'\in G$. If $D\in\mathbb N$
is an upper bound for the diameter
of $K$, and if $g^{-1}g'=kak'$ in the $KAK$-decomposition, then
$d(g,g')-2D\leq d(\mathbf 1,a)\leq d(g,g')+2D$. Since we are only
interested in quasi-isometries, the constant $D$ can be disregarded.
In a similar way, this also yields an $\R$-metric on the ultrapower $G(\R)$ of $G(\bbf R)$ (note that $D$ is also an upper bound for the diameter of
$K(\R)$). We can now explicitly calculate the cone and ultracone from this data:

\begin{Thm}\label{liecone} {\rm (\cite{KT,Thornton})}
Let $G=G(\bbf R)$ be a semisimple  Lie group of non-compact type, and $\rho=e^{-\alpha}$.

Then $Cone(G)=G(^\rho \bbf R)/G(O)$
where $^\rho \bbf R$ is Robinson's valued field $M_0/M_1$ with 
$M_0=\{t\in {\R }|\ |t|<\rho ^{-k} \mbox{ for some } k\in \bbf N\}$ and
$M_1=\{t\in {\R}|\ |t|<\rho ^k \mbox{ for all } k\in \bbf N \}$
and $O=\{\beta\in {^\rho\bbf R}|\ std(log_\rho|\beta|)\geq 0\}$ is its
natural valuation ring.

Also, $UCone(G)=G(\R)/G({\R}_\Fin)$ 
where $\R_\Fin=\{t\in {\R} |\ |t|\leq k \mbox{ for some } k\in \bbf N\}$.
\end{Thm}
The fields $^\rho \bbf R$ were first introduced by Robinson, see \cite{Robinson} and were further studied among others by Luxemburg and Pestov. They will be crucial in Section~\ref{finpres}.

So far, we have just obtained an explicit description of the cones, but in order to be able to use this description, we should obtain more information about its structure. This is given in the following section.

\section{Affine $\Lambda$-buildings and spherical buildings}

Buildings were introduced by Tits to give a (uniform) geometric interpretation of algebraic groups and Lie groups. There are several approaches to introducing them, either as simplicial complexes or as systems of apartments, which will be our point of view in this section.

There are two important major classes of buildings, the spherical buildings, in which the apartments are realized as spheres, and the affine buildings, in which the apartments are realized as affine spaces.

The typical examples of spherical buildings are the ones coming from the standard BN-pair of an algebraic group where the apartments correspond
to the maximal split tori.
In fact, Tits proved that any spherical building of rank at least three arises in this way.
He also proved the following in \cite{TitsLNM} 5.8:

\begin{Thm}\label{Tits} Let $G,G'$ be adjoint absolutely simple
algebraic groups defined over fields $k,k'$, of rank at least $2$.
If the spherical buildings $\Delta$ and $\Delta'$ corresponding to the
groups $G(k)$ and $G'(k')$ are isomorphic, 
then the fields $k\cong k'$ are isomorphic, and $G\cong G'$ are isomorphic as algebraic groups.
\end{Thm}
For affine buildings, the basic example to keep in mind is an infinite tree. There, any branch through the tree can be considered as a copy of the real line.
In a $\Lambda$-tree, the branches of the tree are modeled on a divisible ordered
abelian group $\Lambda$ \cite{LTree}. 
For our purposes, we need to work with a larger class of structures, namely the affine $\Lambda$-buildings, which were introduced by Bennett. These
buildings generalize both affine buildings and $\Lambda$-trees.
In these buildings, the apartments are affine $\Lambda$-spaces 
$\Lambda^n$, for some divisible ordered abelian group $\Lambda$.
The way apartments are required to intersect, is described in the following definition:

\medskip\noindent
\textbf{Affine $\Lambda$-buildings }

\noindent
Let $\Lambda$ be an ordered divisible abelian group, and let  $W$ be a finite reflection group arising from a (crystallographic)
root system. Then $W$ acts naturally on a  lattice $L\cong\bbf Z^n$.

Thus, $W$ acts on $L_\Lambda=L\otimes\Lambda$ as a reflection group,
and the semi-direct product $W_\Lambda=W\ltimes L_\Lambda$
is an affine reflection group acting on the $\Lambda$-metric space $L_\Lambda$ by isometries. A hyperplane fixed by a reflection of $W_\Lambda$ is called a {\em wall}; a wall divides $L_\Lambda$ into two \emph{halfspaces}, and there
is a natural notion of \emph{Weyl chambers}, which are
certain intersections of halfspaces.
\begin{center}
\begin{pspicture}(8,2)
\psframe[linestyle=none,fillstyle=solid,fillcolor=lightgray](3.3,1)(4.7,1.7)
\pswedge[linestyle=none,fillstyle=solid,fillcolor=lightgray](7,1){1}{0}{45}
\psline[linestyle=dotted](.3,.3)(1.7,1.7)
\psline[linewidth=0.8mm](.3,1)(1.7,1)
\psline[linestyle=dotted](1,.3)(1,1.7)
\psline[linestyle=dotted](.3,1.7)(1.7,.3)
\psline[linestyle=dotted](3.3,.3)(4.7,1.7)
\psline[linestyle=dotted](3.3,1)(4.7,1)
\psline[linestyle=dotted](4,.3)(4,1.7)
\psline[linestyle=dotted](3.3,1.7)(4.7,.3)
\psline[linestyle=dotted](6.3,.3)(7.7,1.7)
\psline[linestyle=dotted](6.3,1)(7.7,1)
\psline[linestyle=dotted](7,.3)(7,1.7)
\psline[linestyle=dotted](6.3,1.7)(7.7,.3)
\rput(1,0){a wall}
\rput(4,0){a halfspace}
\rput(7,0){a Weyl}
\rput(7,-.4){chamber}
\end{pspicture}\\\
\end{center}

\begin{Def}{\rm \cite{Bennett}}\label{affinebuilding} The pair $(\mathcal I,\mathcal A)$ is called an 
\emph{affine $\Lambda$-building} of dimension $n$ if 
$\mathcal I$ is a $\Lambda$-metric space and if
the \emph{atlas} $\mathcal A$ on $\mathcal I$ is a family of  $\Lambda$-isometric injections 
$L_\Lambda\too\mathcal I$, called \emph{coordinate charts}, with the following properties:

\textbf{(A1)} If $\phi$ is in $\mathcal A$ and $w\in W_\Lambda$, then
$\phi\circ w:L_\Lambda\too \mathcal I$ is in $\mathcal A$.

\textbf{(A2)} Given two charts $\phi_1,\phi_2$, the set
$B=\phi^{-1}_2(\phi_1(L_\Lambda))$
is convex in $L_\Lambda$, and there exists a $w\in W_\Lambda$
with $\phi_2|_B=\phi_1\circ w|_B$.

The sets $F=\phi(L_\Lambda)\subseteq\mathcal I$ are called \emph{apartments};
the image $S=\phi(S_0)$
of a Weyl chamber $S_0\subseteq L_\Lambda$ is called a \emph{sector}.

\textbf{(A3)} Given $x,y\in\mathcal I$, there exists an apartment
$F=\phi(L_\Lambda)$ containing $x$ and $y$.

\textbf{(A4)} Given two sectors $S_1,S_2\subseteq\mathcal I$, there exist
subsectors $S_1'\subseteq S_1$ and $S_2'\subseteq S_2$ and an apartment
$F$ containing $S_1'\cup S_2'$.

\textbf{(A5)} If $F_1,F_2,F_3$ are apartments such that each of the three
sets $F_i\cap F_j$, $i\neq j$, is a half-apartment
(i.e. the $\phi$-image of a halfspace),
then $F_1\cap F_2\cap F_3\neq\emptyset$.

\textbf{(A6)}
For any apartment $F$ and any $x\in F$, there exists a retraction
$\rho_{x,F}:\mathcal I\too F$ (i.e. $\rho_{x,F}^2=\rho_{x,F}$, and
$\rho_{x,F}$ fixes $F$ pointwise) 
which diminishes distances, with
$\rho_{x,F}^{-1}(x)=\{x\}$.
\end{Def}
Parreau showed that for $\Lambda=\bbf R$ these axioms are equivalent to 
the axioms for a euclidean building given by Kleiner and Leeb in \cite{KL}. The axioms in \cite{KL} are given in terms of geodesics on the metric structure of the underlying set $\mathcal I$ and more differential geometric in flavor. Therefore, we will use the term 'euclidean building' when we want to emphasize that this other set of axioms has been verified, even though this notion is equivalent to an affine $\bbf R$-building.

Bennett also showed that to each affine $\Lambda$-building there is an
associated spherical building, the building at infinity. For $\Lambda=\bbf R$, this building can be visualized in the following way: each affine apartment is a real affine space, say $\bbf R^n$ which can be thought of as having an $(n-1)$-sphere sitting at infinity, the points on the sphere being determined by geodesic rays inside the affine space emanating from some fixed (arbitrary) point. In the case of a tree, the building at infinity is just the set
of ends.

For  affine $\Lambda$-buildings in general there is a combinatorial procedure to associate to each $\Lambda$-apartment of dimension $n$ a combinatorial $(n-1)$-sphere similar to the original way this was done by Tits in \cite{Tits}. Bennett showed that these spheres then form a spherical building, the spherical building at infinity associated to this affine $\Lambda$-building.

So far, very few examples were known except for the
ones coming from the group $SL_n$ over an arbitrary field with valuation. The following result shows that there are indeed many natural  examples of such buildings around.

\begin{Thm}\label{affexamples} (\cite{KT}) Let $R$ be a real closed field, $O\subset R$ an 
o-convex valuation ring (i.e., $O$ is convex and $a\not\in O \Rightarrow a^{-1}\in O$).
Let $G$ be a semisimple Lie group and $A\leq G$ its maximal
(connected) $R$-split torus.

Then $(\mathcal I, \mathcal A)$ is an affine $ \Lambda$-building for $\Lambda =R^*/O^*$, where $\mathcal I=G(R)/G(O)$ and $\mathcal A=\{g:A(R)/A(O)\hookrightarrow \mathcal I|\ g\in G(R)\}$.

The spherical building at infinity of $(\mathcal I, \mathcal A)$ is the standard one coming from the canonical BN-pair of $G(R)$.
\end{Thm}
This result applies in particular in the situation of Theorem \ref{liecone}.
>From the saturation properties of ultrapowers it follows easily that the asymptotic cone of $G$ is a complete metric space and that the atlas $\mathcal A$ one obtains is maximal (i.e., any subset of $Cone(G)$ isometric to $\bbf R^n$ is already in the image of some element of $\mathcal A$).

\begin{Cor}
$Cone(G)$ is a complete affine $\bbf R$-building and 
$UCone(G)$ is an affine $\Lambda$-building for $\Lambda= \R/\R_\Fin)$.
\end{Cor}
In \cite{KL}, Kleiner and Leeb verify that $Cone(G)$ satisfies their axioms of a euclidean building without obtaining an explicit presentation of the building. Hence, using \cite{Parreau} the first part of the corollary already follows from their work. This was crucial in their proof of the Margulis' Conjecture.

Affine $\bbf R$-buildings are rigid in the sense that any homeomorphism preserves the building structure:

\begin{Prop}{\rm \cite{KL,KT}}
\label{TopologicalRigidity}
Let $(\mathcal I,\mathcal A)$ and $(\mathcal I',\mathcal A')$ be complete
affine $\bbf R$-buildings  with maximal atlas $\mathcal A$ and $\mathcal A'$, respectively. Then any  homeomorphism of $\mathcal I$ onto $\mathcal I'$  preserves apartments; moreover, it also induces an isomorphism between the spherical buildings at infinity.
\end{Prop}
Thus we can state the following crucial observation:
\begin{Cor} {\label{crucial}Let $G$ and $G'$ be absolutely simple Lie groups 
of rank at least $2$. The affine $\bbf R$-buildings 
$G(R)/G(O)$ and $ G'(R')/G'(O')$ obtained in \ref{liecone}
are homeomorphic if and only if $G\cong G'$ and $R\cong R'$ as valued fields with valuation rings $O$ and $O'$, respectively.}
\end{Cor}

\section{Application to Margulis' conjecture}

Let  $G=G(\bbf R)$ be a semisimple Lie group of rank $n$ and without compact factors. If $K=K(\bbf R)$ is a maximal compact subgroup, then $X=G/K$ is a Riemannian symmetric space of rank $n$. Every non-compact Riemannian symmetric space is of this form and we may assume that the Riemannian metric on $X$ is that induced by the metric on $G$ defined above. 
For these facts, see \cite{Eberlein}, \cite{Helgason};
it is a well-known fact that
a semisimple Lie group is (essentially) the same as the group of
real points of a semisimple algebraic group defined over $\mathbb Q$,
see \cite{Eberlein} 1.14.
Clearly, $Cone(X)=Cone(G)$, and so we see that the asymptotic cone of a Riemannian symmetric space is an affine $\bbf R$-building. 
Under the cone construction, the {\em maximal flats} of $X$, i.e. the subspaces of $X$ isometric to $\bbf R ^n$ are transformed into the apartments of the affine building.

Margulis conjectured that if $X$ and $Y$ are any two Riemannian symmetric spaces without compact factors and $f:X\too Y$ is a quasi-isometry, then there is an isometry $\tilde f: X\too Y$ at bounded distance from $f$.

The conjecture was proved for Riemannian symmetric spaces $X=G/K$ where $G$ contains no simple factors of rank at most $1$ by Kleiner and Leeb. (There are counterexamples in rank $1$, and the situation was completely described by Pansu in \cite{Pansu}.) They used the fact, that the asymptotic cone of $X$ is a euclidean building and work with the geodesic structure of $Cone(X)$.

\medskip\noindent\textbf{Proof of the Margulis conjecture}

\noindent
We can now outline a proof of the Margulis conjecture, using the previous
results. The details will appear in \cite{KT}.
Let $X=G/K$ and $X'=G'/K'$ be Riemannian symmetric spaces of noncompact type,
where $G$ and $G'$ contain no factors of rank at most $1$, and let 
\[
f:X\too X'
\]
be an $(L,C)$-quasi-isometry.
Using the definition of the metric it can be shown that we may interpret the metric spaces $(X,d)$ and $(X,d')$ in the structure $(\bbf R,+,\cdot,\log)$. Thus, we may take an ultrapower of this structure extended by the quasi-isometry $f$ to obtain
an $(L,C)$ quasi-isometry
\[
^*f:{}^*X\too{}^*X'
\]
between the $\R$-metric spaces $^*X,{}^*X'$. 

We then obtain a (bi-Lipschitz) homeomorphism
\[
Cone(f):Cone(X)=G(R)/G(O)\too G'(R)/G'(O)=Cone(X')
\]
on the asymptotic cones, where $R={}^\rho\bbf R$
denotes Robinson's field, and
$O$ its natural valuation ring.

By Theorem \ref{TopologicalRigidity} we know that $Cone(f)$ also induces an isomorphism between the spherical buildings at infinity of $Cone(X)$ and $Cone(X')$.
Using Theorem \ref{Tits} we
can now already conclude that the groups $G$ and $G'$ are isomorphic, which implies that $X$ and $X'$ are essentially the same. 

In order to find an isometry at bounded distance from $f$, we use the ultracone $UCone(X)$ of $X$.

First we show that there is some finite bound $c\in\mathbb N$
such that for any maximal flat $A$ of $^*X$ (i.e. a subset $A$ of $^*X$ which is $\R$-isometric to $\R^n$) there is a flat $A'$ such that the image of $A$ under $^*f$ has Hausdorff distance at most $c$ from $A'$. The main point here is that any two apartments
in an affine building having finite Hausdorff distance are equal, and
that $Cone(f)$ takes apartments to apartments for \emph{any}
asymptotic cone.
Such an 'approximating' flat $A'$ exists thus on every large metric ball in
${}^*X'$, and, by repeated application of compactness arguments,
we deduce the global existence of the flat $A'$. In particular, $^*f$ induces
a well-defined bijection between the flats of ${}^*X$ and ${}^*X'$.
Since the maximal flats in $X$ and $^*X$ can be described using the group structure of $G(\bbf R)$ and $G(\R)$, respectively, this map 
is first order expressible and exists also for the
quasi-isometry $f:X\rTo X'$.

Since $c$ was a finite constant, $UCone(f)$ is an apartment preserving map.
By \ref{affexamples}, $UCone(X)$ is the affine $\R/\R_\Fin$-building $G(\R)/G(\R_\Fin)$ whose spherical building at infinity is that of $G(\R)$. Using nonstandard analysis and the structure of this affine  $\R/\R_\Fin$-building we can show that $UCone(f)$ induces a continuous map on the spherical buildings at infinity of $UCone(X)$ and $UCone(X')$
(which are the canonical spherical buildings associated to $G(\R)$ and
$G'(\R)$). Using the map defined
on the flats, it is not difficult to see that this continuous map
is also first order expressible; thus, the quasi-isometry $f$
induces a continuous map on the spherical buildings associated to
$G(\bbf R)$ and $G'(\mathbb R)$. For rank at least $2$, the group of continuous automorphisms of this spherical building coincides with the isometry group of $X$ (and is in fact a finite extension of  $G$). This is essentially due
to the fact that every continuous automorphism of a semisimple Lie group
is algebraic, and that the group can be recovered from the building.

Thus, there is some 
isometry $g:X\rTo X'$ inducing on the spherical buildings of $G(\bbf R)$ 
and $G'(\mathbb R)$ the same map as the quasi-isometry $f$. 
It follows that $UCone(f)$ and $UCone(g)$ agree (set-wise) on the
collection of all apartments
of $UCone(X)$. But this implies that $UCone(f)=UCone(g)$. By Remark~\ref{remark}, $f$ and $g$ have bounded distance and we are done.

\section{\label{finpres} Application to finitely presented groups}

We now turn to Gromov's question  whether there exists a finitely presented group with non-homeomorphic asymptotic cones. This is closely related to asymptotic cones of Lie groups via the following characterization of finitely presented groups (see \cite{BH}  p.137):

\begin{Thm}{ A group $\Gamma$ is finitely presented if and only if it acts properly and cocompactly by isometries on a simply connected geodesic space $X$.}
\end{Thm}
Here, an action is called {\em cocompact} if there exists a compact set $K\subset X$ such that $X=\Gamma .K$. The action is called {\em proper} if there is some number $r>0$ such that for each $x\in X$ the set $\{\gamma\in\Gamma\ |\ \gamma(B(x,r))\cap B(x,r)\neq\emptyset\}$ is finite.

Riemannian symmetric spaces of noncompact type
are in particular simply connected and geo\-desic, and for these spaces subgroups of the isometry group which act cocompactly are also called {\em uniform lattices}. By an old result of Borel, any Riemannian symmetric space allows many such uniform lattices. Note that such a lattice is quasi-isometric to the Riemannian symmetric space on which it acts and hence its asymptotic cone is homeomorphic to  that of the symmetric space. Thus, if we can produce non-homeomorphic asymptotic cones of Riemannian symmetric spaces, we have also produced examples of finitely presented groups with non-homeomorphic asymptotic cones.

Since we know that the asymptotic cone of the Riemannian symmetric space $X$ is the affine $\bbf R$-building
$G(R)/G(O)$,
with $R$ the Robinson field ${}^\rho\bbf R$
and $O$ its valuation ring, it suffices by Proposition~\ref{crucial} to produce ultrafilters giving rise to non-isomorphic Robinson fields.
It is not hard to see that if we assume \textsf{CH}, then all Robinson fields are isomorphic since any countable ultrapower is saturated. Hence in this case it is impossible to produce non-homeomorphic asymptotic cones from finitely presented groups acting on Riemannian symmetric spaces. 

So we now assume $\neg$\textsf{CH}.
Examples constructed by Roitman \cite{roitman} yield non-isomorphic ultrapowers differing from each other through their degree of saturation and the existence of scales. However, in the Robinson field we truncate behind some large non-standard number $\alpha$ and divide out by the relative 'infinitesimal' part, which is exactly the part which contains the saturation.

Hence, in order to produce non-isomorphic asymptotic cones we concentrate on producing ultrafilters which give rise to ultrapowers of the natural numbers which differ in the coinitial segments of $\omega$.

It is not hard to see that this suffices in order to make sure that the Robinson fields obtained from these ultrafilters will be non-isomorphic.

Thus, assuming $\neg$\textsf{CH}, the following result shows that there are finitely presented groups with non-homeomorphic asymptotic cones:

\begin{Thm}{{\rm \cite{KSTT}} Assume $\neg$\textsf{CH}. Then there are $2^{2^{\omega}}$  Robinson fields up to isomorphism.}
\end{Thm}

In fact, we show the following:
\begin{Cor}{{\rm \cite{KSTT}} Let $G$ be a connected semisimple Lie group with at least
one absolutely simple factor $S$ such that 
$\bbf R$-${\rm rank}(S) \geq 2$ and let $\Gamma$ be a
uniform lattice in $G$.
\begin{enumerate}
\item[(a)] If \textsf{CH} holds, then $\Gamma$ has a unique
asymptotic cone up to homeomorphism.
\item[(b)] If \textsf{CH} fails, then $\Gamma$ has 
$2^{2^{\omega}}$ asymptotic cones up to homeomorphism.
\end{enumerate}}
\end{Cor}

Furthermore,
\begin{Thm}{{\rm \cite{KSTT}} 
If \textsf{CH} holds, then every finitely generated group $\Gamma$ 
has at most $2^{\omega}$ asymptotic cones up to isometry.}
\end{Thm}

\begin{Cor} {{\rm \cite{KSTT}} 
The following statements are equivalent.
\begin{enumerate}
\item[(a)] \textsf{CH} fails.
\item[(b)] There exists a finitely presented group $\Gamma$
which has $2^{2^{\omega}}$ asymptotic cones up to
homeomorphism.
\item[(c)] There exists a finitely generated group $\Gamma$
which has $2^{2^{\omega}}$ asymptotic cones up to
homeomorphism.
\end{enumerate}}
\end{Cor}

\small\sc\noindent 
Linus Kramer,
Fachbereich Mathematik,
TU Darmstadt,
Schlo\ss gartenstra\ss e~7,
D--64289 Darmstadt,
Germany
{\tt kramer@mathematik.tu-darmstadt.de}\\ 
\\
Katrin Tent,
Mathematisches Institut,
Universit\"at W\"urzburg,
Am Hubland,
D--97074 W\"urzburg,
Germany
{\tt tent@mathematik.uni-wuerzburg.de} 
 
\end{document}